\documentclass{amsart}

\usepackage{graphicx}

\theoremstyle{definition}

\theoremstyle{remark}

\numberwithin{equation}{section}

\begin{document}

\title{The chromatic number of the plane is at least 5}
\author{Aubrey D.N.J. de Grey}

\address{SENS Research Foundation, Mountain View, California 94041, USA; email: aubrey@sens.org}

\begin{abstract}
We present a family of finite unit-distance graphs in the plane that are not 4-colourable, thereby improving the lower bound of the Hadwiger-Nelson problem. The smallest such graph that we have so far discovered has 1581 vertices.
\end{abstract}

\maketitle

\section{Introduction}
How many colours are needed to colour the plane so that no two points at distance exactly 1 from each other are the same colour? This quantity, termed the chromatic number of the plane or CNP, was first discussed (though not in print) by Nelson in 1950 (see [Soi]). Since that year it has been known that at least four and at most seven colours are needed. The lower bound was also noted by Nelson (see [Soi]) and arises because there exist 4-chromatic finite graphs that can be drawn in the plane with each edge being a straight line of unit length, the smallest of which is the 7-vertex Moser spindle [MM] (see Figure 7, left panel). The upper bound arises because, as first observed by Isbell also in 1950 (see [Soi]), congruent regular hexagons tiling the plane can be assigned seven colours in a pattern that separates all same-coloured pairs of tiles by more than their diameter. The question of the chromatic number of the plane is termed the Hadwiger-Nelson problem, because of the contributions of Nelson just mentioned and because the 7-colouring of the hexagonal tiling was first discussed (though in another context) by Hadwiger in 1945 [Had]. The rich history of this problem and related ones is wonderfully documented in [Soi]. Since 1950, no improvement has been made to either bound.

We now give a brief overview of our construction. The graphs described in this report are obtained as follows:
\begin{enumerate}
\item	We note that the 7-vertex, 12-edge unit-distance graph $H$ consisting of the centre and vertices of a regular hexagon of side-length 1 can be coloured with at most four colours in four essentally distinct ways (that is, up to rotation, reflection and colour transposition). Two of these colourings contain a monochromatic triple of vertices and two do not.
\item	We construct a unit-distance graph $L$ that contains 52 copies of $H$ and show that, in all 4-colourings of $L$, at least one copy of $H$ contains a monochromatic triple.
\item	We construct a unit-distance graph $M$ that contains a copy of $H$ and show that there is no 4-colouring of $M$ in which that $H$ contains a monochromatic triple. Thus, the unit-distance graph $N$ created by arranging 52 copies of $M$ so that their counterparts of $H$ form a copy of $L$ is not 4-colourable. This completes our demonstration that the CNP is at least 5.
\item	We identify smaller non-4-colourable unit-distance graphs, first by identifying vertices in $N$ whose deletion preserves the absence of a 4-colouring, and then by more elaborate methods.
\end{enumerate}

\section{The 4-colourings of $H$}
Figure 1 shows the four essentially distinct (i.e., up to rotation, reflection and colour transposition) colourings of $H$ with at most four colours. The top two colourings possess, and the bottom two lack, a triple of vertices all the same colour.

\begin{figure}
\includegraphics[width=4in]{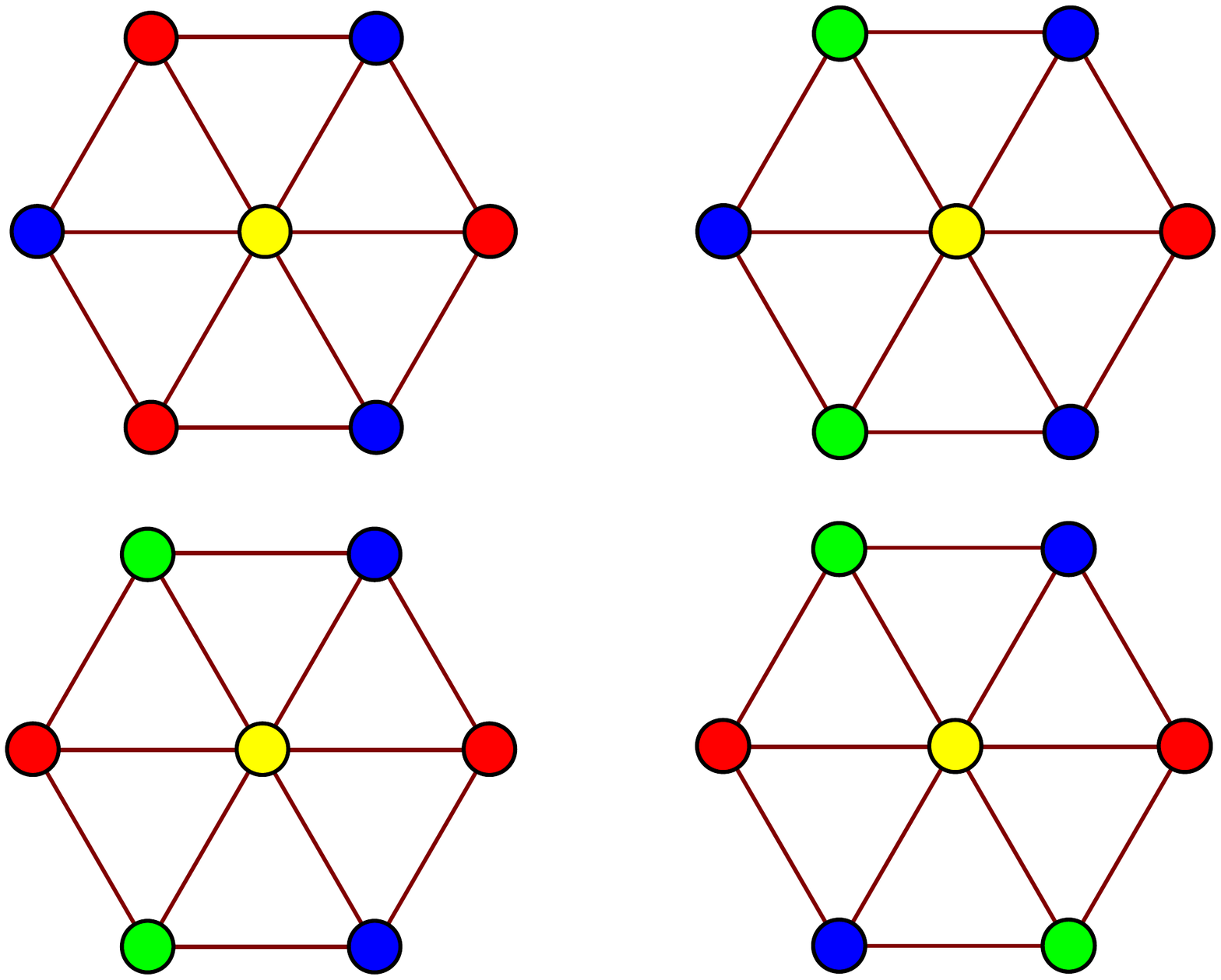}
\caption{The essentially distinct ways to colour $H$ with at most four colours.}
\label{firstfig}
\end{figure}

\section{Construction and colourings of $L$}
\subsection{A 31-vertex graph $J$ assembled from 13 copies of $H$}

The graph $J$, shown in Figure 2, contains a copy of $H$ in its centre, six copies centred at distance 1 from its centre, and six copies centred at distance \(\sqrt{3}\) from its centre.

\begin{figure}
\includegraphics[width=5in]{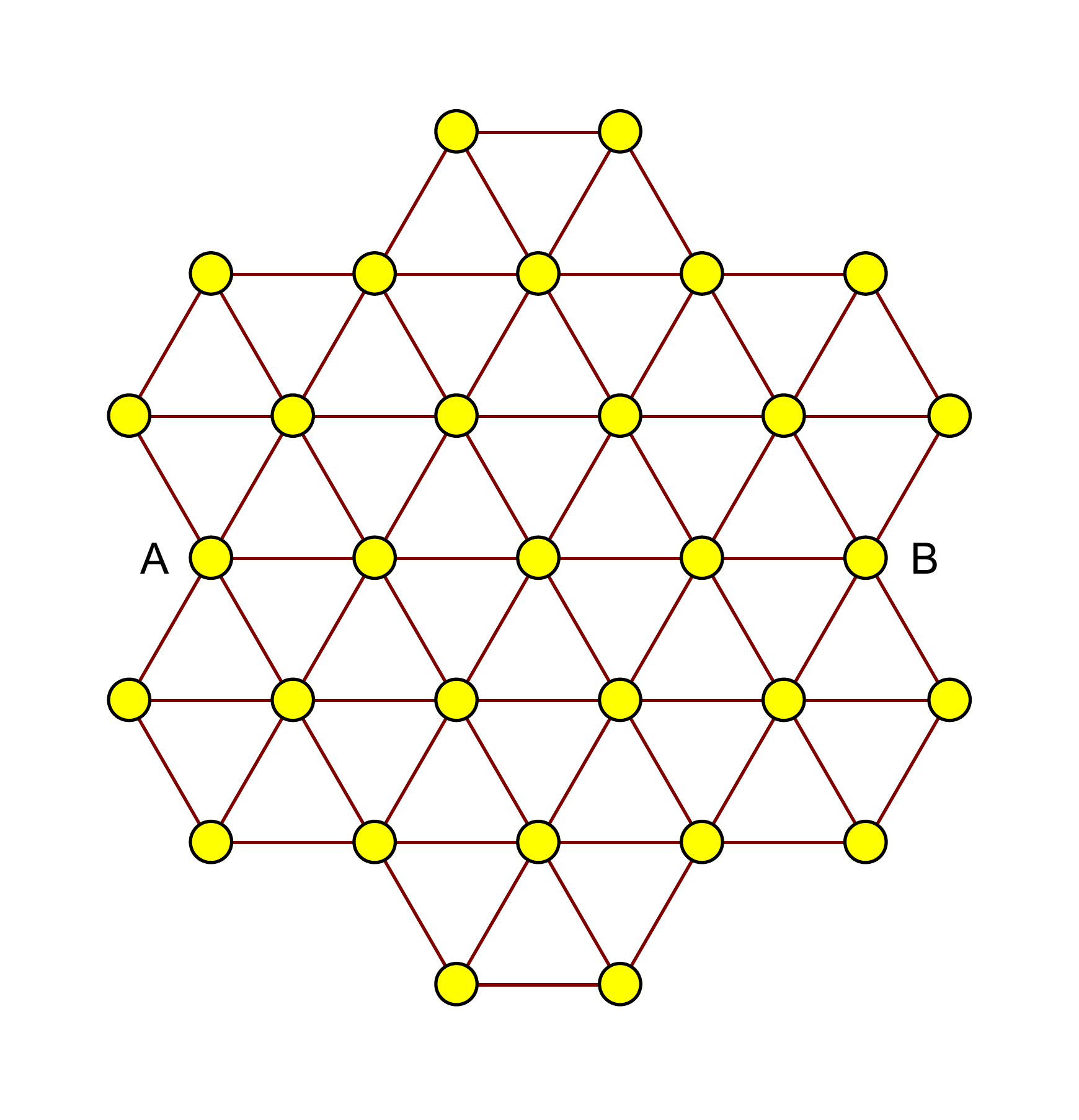}
\caption{The graph $J$, containing 31 vertices and 13 copies of $H$.}
\label{secondfig}
\end{figure}

\subsection{The 4-colourings of $J$ in which no copy of $H$ contains a monochromatic triple}

Figure 3 shows six essentially distinct such colourings up to rotation, reflection, colour transposition and the colours of the vertices coloured black in Figure 3 (which will not concern us hereafter). That these are the only such colourings can be checked by grouping the possibilities according to whether the central copy of $H$ has two (the top row in Figure 3) or no (bottom row) monochromatic pairs of vertices at distance \(\sqrt{3}\) and, in the case where it has none, whether all the copies of $H$ centred at distance 1 from the centre also have none (bottom left colouring) or some have two (last two colourings). The bottom centre colouring is the reason we need the vertices coloured black; though those vertices can be coloured in many ways, it turns out that if they were deleted then there would be additional 4-colourings of the remaining graph in which none of the seven remaining copies of $H$ contained a monochromatic triple, and some of those new colourings lack a key property shared by all those in Figure 3, which we now discuss.

\begin{figure}
\includegraphics[width=5in]{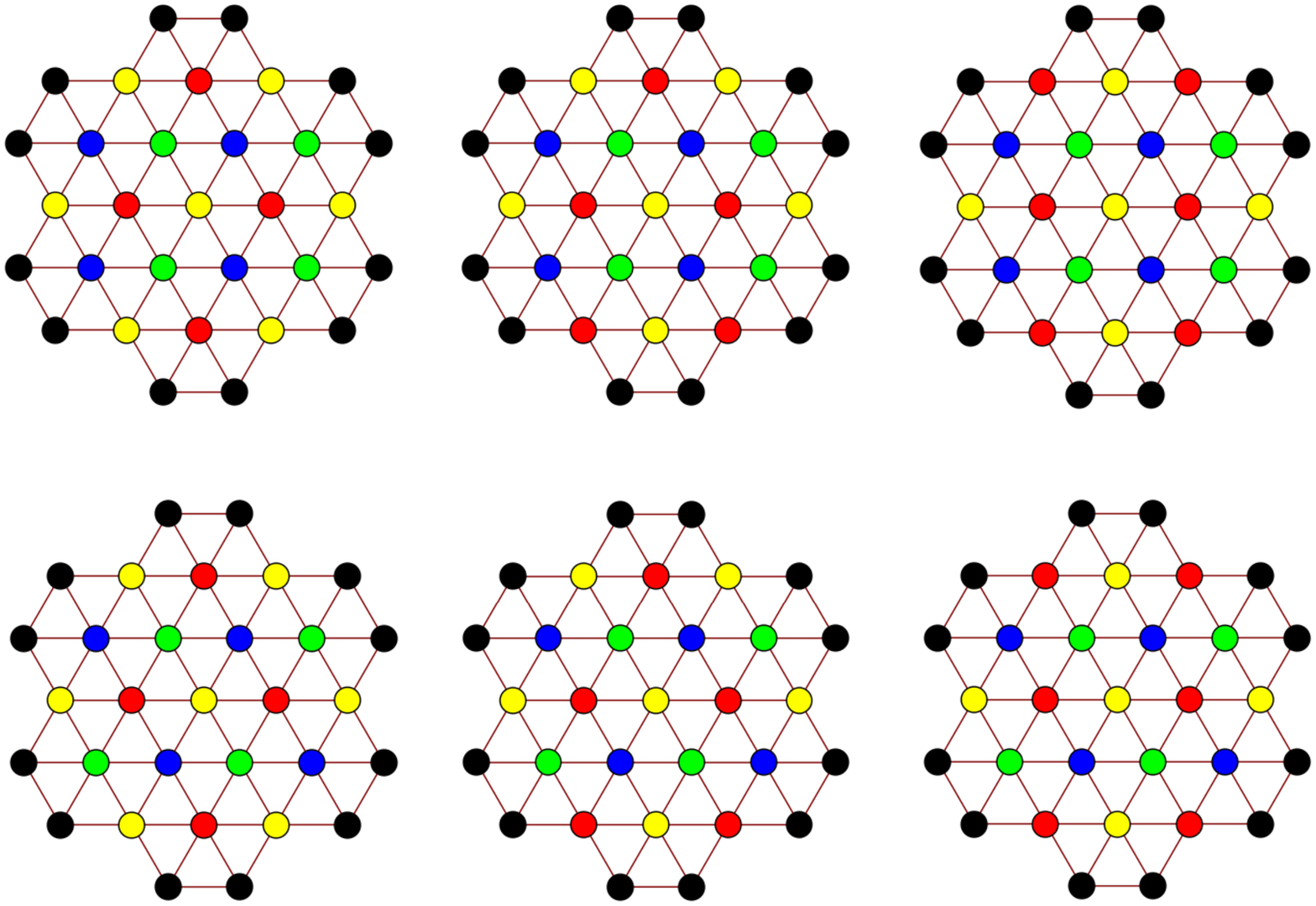}
\caption{The essentially distinct colourings of $J$ in which no copy of $H$ (including the ones containing some black vertices) contains a monochromatic triple. Black vertices can assume any colour so long as no connected vertices are the same colour.}
\label{thirdfig}
\end{figure}

The aspect of the colourings in Figure 3 that will be our focus below is that they feature only three essentially distinct colourings of the vertices at distance 2 from the centre, which we shall hereafter term the linking vertices. Either (a) the linking vertices are all the same colour as the centre (left-hand cases in Figure 3), or (b) four consecutive ones (when enumerated going clockwise around the centre) are the same colour as the centre and the other two are a second colour (middle cases), or (c) two opposite ones are the same colour as the centre and all the other four are a second colour (right-hand cases). Hereafter we refer to a pair of linking vertices located in opposite directions from the centre, such as those labelled $A$ and $B$ in Figures 2 and 4, as a linking diagonal.

\subsection{A 61-vertex graph $K$ assembled from two copies of $J$}

We now construct $K$ as the union of $J$ with a copy of $J$ rotated clockwise around the origin by \(2 \arcsin(1/4)\). This rotation causes corresponding linking vertices to lie at distance 1 from each other; see Figure 4.

\begin{figure}
\includegraphics[width=5in]{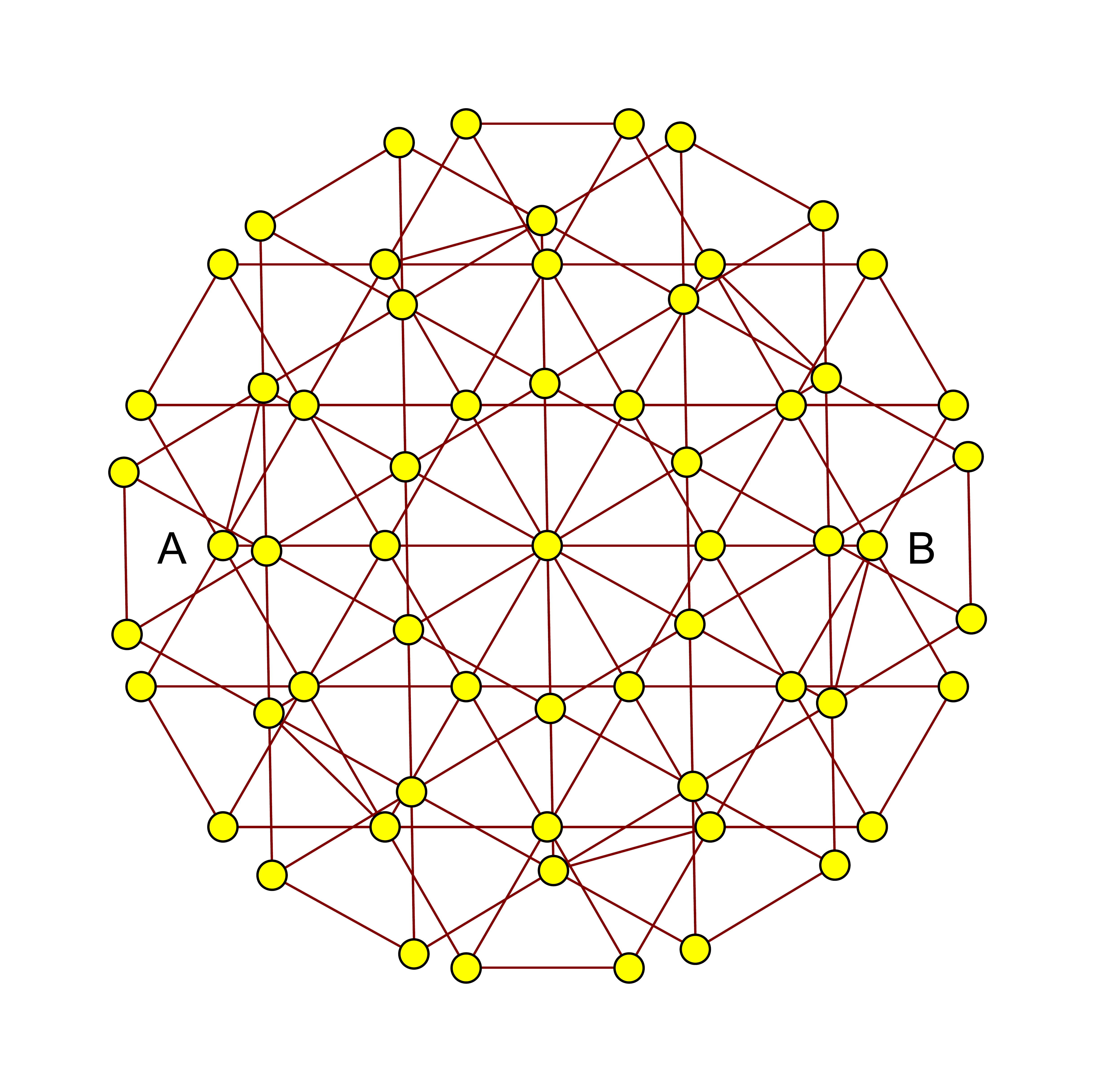}
\caption{The graph $K$, containing 61 vertices and 26 copies of $H$.}
\label{fourthfig}
\end{figure}

We now note that, in any 4-colouring of $K$ in which none of the 26 copies of $H$ contains a monochromatic triple, both copies of $J$ must have their linking vertices coloured according to option (c) above. This is of interest because in option (c) each of the three linking diagonals of $J$ is monochromatic. Thus, in all 4-colourings of $K$ where no copy of $H$ contains a monochromatic triple, all six linking diagonals are monochromatic.

\subsection{A 121-vertex graph $L$ assembled from two copies of $K$}

Finally, we construct $L$ as the union of $K$ with a copy of $K$ rotated around $A$ by \(2 \arcsin(1/8)\). This rotation causes the counterpart of $B$ (denoted $B'$) to lie at distance 1 from $B$ (see Figure 5, in which $L$ has been translated and rotated so as to give it symmetry about the y-axis).

\begin{figure}
\includegraphics[width=5in]{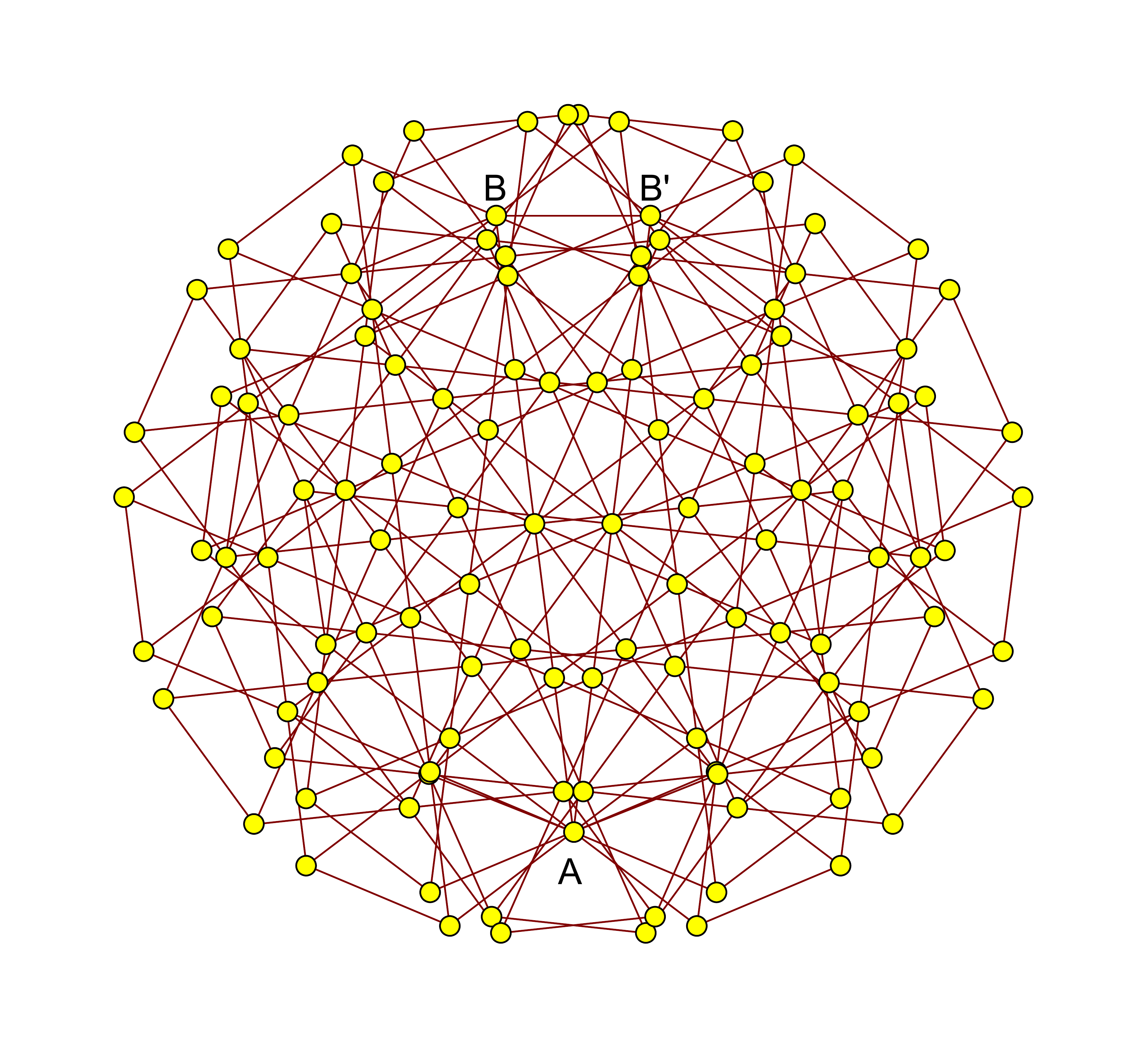}
\caption{The graph $L$, containing 121 vertices and 52 copies of $H$.}
\label{fifthfig}
\end{figure}

The property noted above for $K$ guarantees that in no 4-colouring of $L$ do all of its 52 constituent copies of $H$ lack a monochromatic triple. Either $B$ or $B'$ must be a different colour than $A$, so one of the copies of $K$ must contain a non-monochromatic linking diagonal, so it must contain a copy of $H$ with a monochromatic triple.

\section{Construction and colourings of $M$}

In seeking graphs that can serve as $M$ in our construction, we focus on graphs that contain a high density of Moser spindles. The motivation for exploring such graphs is that a spindle contains two pairs of vertices distance \(\sqrt{3}\) apart, and these pairs cannot both be monochromatic. Intuitively, therefore, a graph containing a high density of interlocking spindles might be constrained to have its monochromatic \(\sqrt{3}\)-apart vertex pairs distributed rather uniformly (in some sense) in any 4-colouring. Since such graphs typically also contain regular hexagons of side-length 1, one might be optimistic that they could contain some such hexagon that does not contain a monochromatic triple in any 4-colouring of the overall graph, since such a triple is always an equilateral triangle of edge \(\sqrt{3}\) and thus constitutes a locally high density, i.e. a departure from the aforementioned uniformity, of monochromatic \(\sqrt{3}\)-apart vertex pairs.

\subsection{Graphs with high edge density and spindle density}

In seeking graphs with high spindle density, we begin by noting two attractive features of the 9-vertex unit-distance graph $T$ that is obtained by adding two particular vertices to a spindle (see Figure 6, left). These vertices, $P$ and $Q$:
\begin{enumerate}
\item form an equilateral triangle with the “tip” $X$ of the spindle
\item lie on (the extension of) the line forming the “base” $YZ$ of the spindle
\end{enumerate}
We can thus construct a 15-vertex unit-distance graph $U$ (Figure 6, right) that contains three spindles and possesses rotational and reflectional triangular symmetry. This symmetry suggests that graphs formed by combining translations and 60-degree rotations of $U$ might have particularly high edge and spindle density.
              
\begin{figure}
\includegraphics[width=4in]{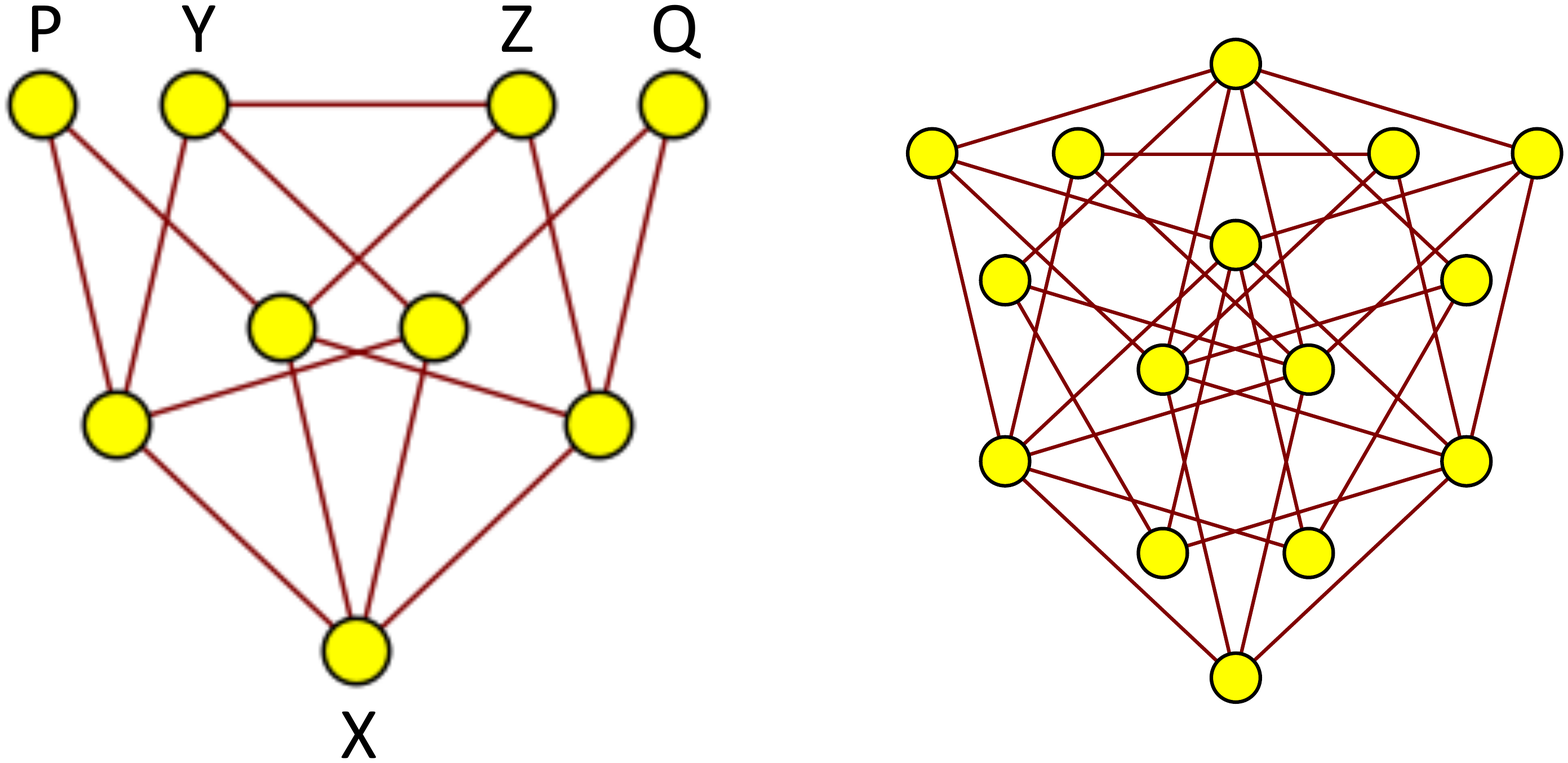}
\caption{The graphs $T$ and $U$.}
\label{sixthfig}
\end{figure}

\subsection{Construction of a graph that serves as $M$}

The speculation just mentioned turns out to be true; for example, we have found a 97-vertex graph containing 78 spindles (not shown). Therefore, a custom program (outlined in the next section) was written to test graphs of this form for possession of a 4-colouring in which the central $H$ contains a monochromatic triple. Disappointingly, we could not find a graph of this form that enforces sufficient uniformity of the distribution of \(\sqrt{3}\)-apart monochromatic vertex pairs to deliver the property we require for $M$, even though we checked examples with well over 1000 vertices.

However, this approach can be extended. Graphs arising from the construction described thus far can have spindles in only six different orientations, with edges falling into just three equivalence classes if equivalence is defined as rotation by a multiple of 60 degrees; see Figure 7 (left). A natural elaboration is to add more classes based on the relative orientations of spindles that share a lot of vertices, such as in Figure 7 (middle). The maximum possible degree of a vertex in a graph constructed from these vectors is thereby increased from 18 to 30; see Figure 7 (right), which we denote as graph $V$. The angles of the edges relative to the vector (1,0) are \(i \arcsin({\sqrt 3}/2) + j \arcsin(1/\sqrt{12}), i \in 0 \dots 5, j \in -2 \dots 2\}\).

\begin{figure}
\includegraphics[width=5in]{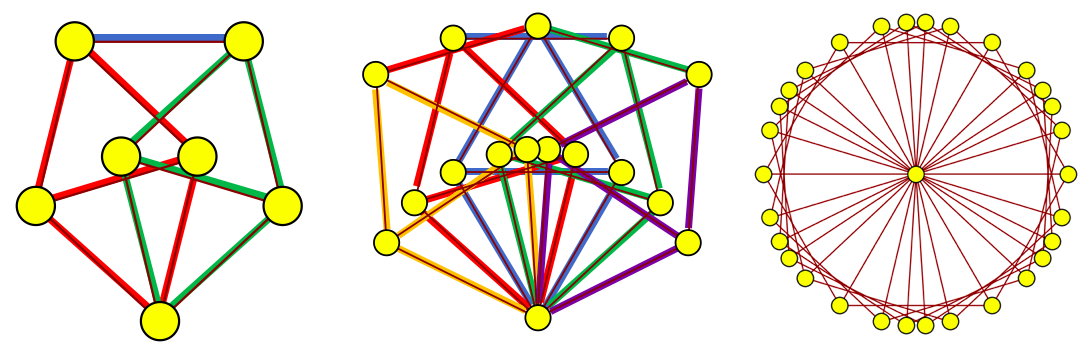}
\caption{The vector classes present in one (left), or three tightly linked (middle), Moser spindles; the graph $V$ (right).}
\label{seventhfig}
\end{figure}

Happily, this turns out to suffice. Let $W$ be the 301-vertex graph consisting of all points at distance \(\leq \sqrt 3\) from the origin that are the sum of two edges of $V$ (interpreted as vectors). The 1345-vertex graph shown in Figure 8 is the union of $W$ with six translates of it in which the origin is mapped to a vertex of $H$. Our program did not find any 4-colouring of this latter graph in which the central $H$ contains a monochromatic triple, so it can serve as our $M$. We can, in other words, create a non-4-colourable unit-distance graph $N$ as the union of 52 copies of $M$, translated and rotated so that each instance of $H$ in $L$ coincides with the central $H$ of a copy of $M$. After merging coincident vertices arising from different copies of $M$, this graph has 20425 vertices. (We do not show a picture of $N$, because it is visually impenetrable and we have anyway discovered smaller examples, discussed below. Even $M$ is too large for its picture to be explanatory; we include it just in case it is interesting to see the general shape.)

\begin{figure}
\includegraphics[width=5in]{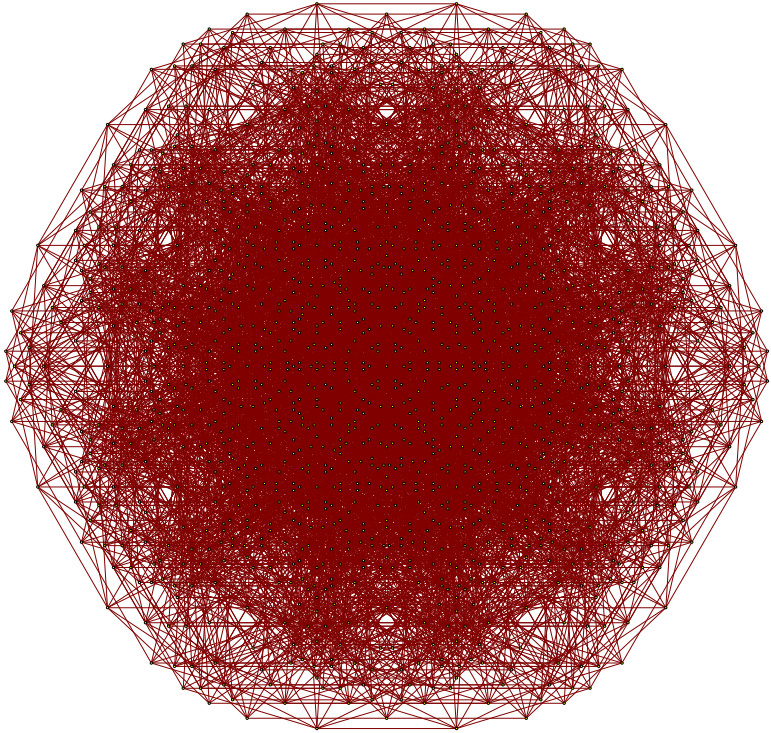}
\caption{The graph $M$.}
\label{eighthfig}
\end{figure}

\subsection{Testing 4-colourability of edge-dense, spindle-dense graphs}

It is, in general, computationally challenging to determine the chromatic number of a graph with over 1000 vertices by simplistic search methods, let alone one with 20425 vertices. We thus developed a custom program to test graphs for possession of the property required for our graph $M$, taking advantage of certain properties of our candidate graphs.

Because we are only asking whether a specific number of colours is or is not sufficient, and also because of the high density of edges and spindles in our target graphs, the required test turns out to be computationally far cheaper than a general determination of chromatic number of comparable-sized graphs. It can be performed rapidly by a simple depth-first search optimised only slightly, as follows:
\begin{enumerate}
\item Since our question is whether there is any 4-colouring of $M$ in which the central $H$ contains a monochromatic triple, we begin by specifying the colourings of the vertices of that central $H$, which we term the initialising vertices. Since $M$ has the same symmetry as $H$, we need only check the two essentially distinct triple-containing colourings of those vertices.
\item We order the remaining vertices according to a hierarchy (most significant first) of parameters (all decreasing): how many spindles they are part of, their degree, and how many unit triangles they are part of.
\item We colour the next not-yet-coloured vertex (initially vertex 8) with the first colour that we have not already tried for it (initially colour 1).
\item We check each not-yet-coloured neighbour (if any) of the just-coloured vertex to see how many colours are still permissible for it. If any such vertex already has neighbours of all four colours, we will need to backtrack (see step 6 below). If any has neighbours of exactly three colours, we assign the remaining (forced) colour to it.
\item If we do not need to backtrack, but we did colour some vertices in step 4, we repeat step 4 for each such vertex.
\item If we need to backtrack, we uncolour everything that we just coloured in the most recent iteration of steps 3 and (any resulting iterations of) 4.
\item If we just did an uncolouring and the just-uncoloured vertex that was coloured at step 3 has no colours that have not yet been tried, we repeat step 6 for the next-most-recent iteration of steps 3 and 4 unless we have already backtracked all the way to the vertices of $H$. Otherwise, if there are still some uncoloured vertices we return to step 3.
\item We terminate when we get here, i.e. when either all vertices are coloured or we have backtracked all the way down to the vertices of $H$.
\end{enumerate}

This algorithm was implemented in Mathematica 11 on a standard MacBook Air and terminated in only a few minutes for our candidate $M$, without finding a 4-colouring starting from either of the triple-containing colourings of its central $H$. Essentially, the speed arises because the fixing of only 20 or so colours at step 3 typically lets almost all remaining colours be forced at step 4.

\section{Identification of smaller solutions}

One may naturally wonder how large the smallest non-4-colourable unit-distance graph is. Moreover, if substantially smaller graphs were found, they might be feasibly verifiable independently by standard algorithms such as SAT solvers, thus avoiding recourse to potentially buggy custom software.

The most direct way to find such graphs is to seek a succession of small simplifications of $N$. Many approaches to this are evident, some much more computationally tractable than others. We have thus far employed only a small range of strategies that take advantage of the stepwise method by which $N$ was constructed: briefly, we identified vertices of subgraphs such as $M$ whose removal preserves the property required of them in the construction and thus also the chromatic number of $N$, and then we sought individual new vertices whose addition allowed the removal of more than one pre-existing vertex.

\subsection{Status of shrinking $N$}

Despite the rather simplistic nature of these methods, we have so far shrunk $N$ by a factor of nearly 13, our current record being the 1581-vertex graph $G$ that (again more for artistic than expository reasons) is shown in Figure 9. It can be constructed as follows:

\begin{enumerate}
\item Let $S$ be the following set of points:

\[(0, 0), (1/3, 0), (1, 0), (2, 0), ((\sqrt{33} - 3)/6, 0),
(1/2, 1/\sqrt{12}), (1, 1/\sqrt{3}), (3/2, \sqrt{3}/2),\]
\[(7/6, \sqrt{11}/6),
(1/6, (\sqrt{12} - \sqrt{11})/6), (5/6, (\sqrt{12} - \sqrt{11})/6),\]
\[(2/3, (\sqrt{11} - \sqrt{3})/6), (2/3, (3\sqrt{3} - \sqrt{11})/6),(\sqrt{33}/6, 1/\sqrt{12}),\]
\[((\sqrt{33} + 3)/6, 1/\sqrt{3}),
((\sqrt{33} + 1)/6, (3\sqrt{3} - \sqrt{11})/6),
((\sqrt{33} - 1)/6, (3\sqrt{3} - \sqrt{11})/6),\]
\[((\sqrt{33} + 1)/6, (\sqrt{11} - \sqrt{3})/6),
((\sqrt{33} - 1)/6, (\sqrt{11} - \sqrt{3})/6),\]
\[((\sqrt{33} - 2)/6, (2\sqrt{3} - \sqrt{11})/6),
((\sqrt{33} - 4)/6, (2\sqrt{3} - \sqrt{11})/6),\]
\[((\sqrt{33} + 13)/12, (\sqrt{11} - \sqrt{3})/12),
((\sqrt{33} + 11)/12, (\sqrt{3} + \sqrt{11})/12),\]
\[((\sqrt{33} + 9)/12, (\sqrt{11} - \sqrt{3})/4),
((\sqrt{33} + 9)/12, (3\sqrt{3} + \sqrt{11})/12),\]
\[((\sqrt{33} + 7)/12, (\sqrt{3} + \sqrt{11})/12),
((\sqrt{33} + 7)/12, (3\sqrt{3} - \sqrt{11})/12),\]
\[((\sqrt{33} + 5)/12, (5\sqrt{3} - \sqrt{11})/12),
((\sqrt{33} + 5)/12, (\sqrt{11} - \sqrt{3})/12),\]
\[((\sqrt{33} + 3)/12, (3\sqrt{11} - 5\sqrt{3})/12),
((\sqrt{33} + 3)/12, (\sqrt{3} + \sqrt{11})/12),\]
\[((\sqrt{33} + 3)/12, (3\sqrt{3} - \sqrt{11})/12),
((\sqrt{33} + 1)/12, (\sqrt{11} - \sqrt{3})/12),\]
\[((\sqrt{33} - 1)/12, (3\sqrt{3} - \sqrt{11})/12),
((\sqrt{33} - 3)/12, (\sqrt{11} - \sqrt{3})/12),\]
\[((15 - \sqrt{33})/12, (\sqrt{11} - \sqrt{3})/4),
((15 - \sqrt{33})/12, (7\sqrt{3} - 3\sqrt{11})/12),\]
\[((13 - \sqrt{33})/12, (3\sqrt{3} - \sqrt{11})/12),
((11 - \sqrt{33})/12, (\sqrt{11} - \sqrt{3})/12)\]

\item Let $S_a$ be the unit-distance graph whose vertices consist of all points obtained by rotating the points in $S$ around the origin by multiples of 60 degrees and/or by negating their $y$-coordinates. $S_a$ has a total of 397 vertices.

\item Let $S_b$ be $S_a$ rotated anticlockwise about the origin by \(2 \arcsin(1/4)\).

\item Let $Y$ be the union of $S_a$ and $S_b$ with the vertices \((1/3,0)\) and \((-1/3,0)\) deleted.

\item Rotate $Y$ anticlockwise about (-2,0) by \(\pi/2 + \arcsin(1/8)\) to give $Y_a$.

\item Rotate $Y$ anticlockwise about (-2,0) by \(\pi/2 - \arcsin(1/8)\) to give $Y_b$.

\item Let $G$ be the union of $Y_a$ and $Y_b$.
\end{enumerate}

Happily, $G$ has turned out to be within the reach of standard SAT solvers, with which others have now confirmed its chromatic number to be 5 without the need to resort to using custom code or checking weaker properties of subgraphs.

This attempt to identify smaller examples has clearly been rather cursory, so we think it highly likely that examples smaller than $G$ exist. Indeed, a Polymath project [DHJP] has been created to seek such graphs, as well as to seek ones whose lack of a 4-colouring can be shown without a computer. Concise and explicit descriptions of certain 5-chromatic unit-distance graphs would seem to be a promising way to attack the question of whether 6-chromatic examples exist.

\begin{figure}
\includegraphics[width=5in]{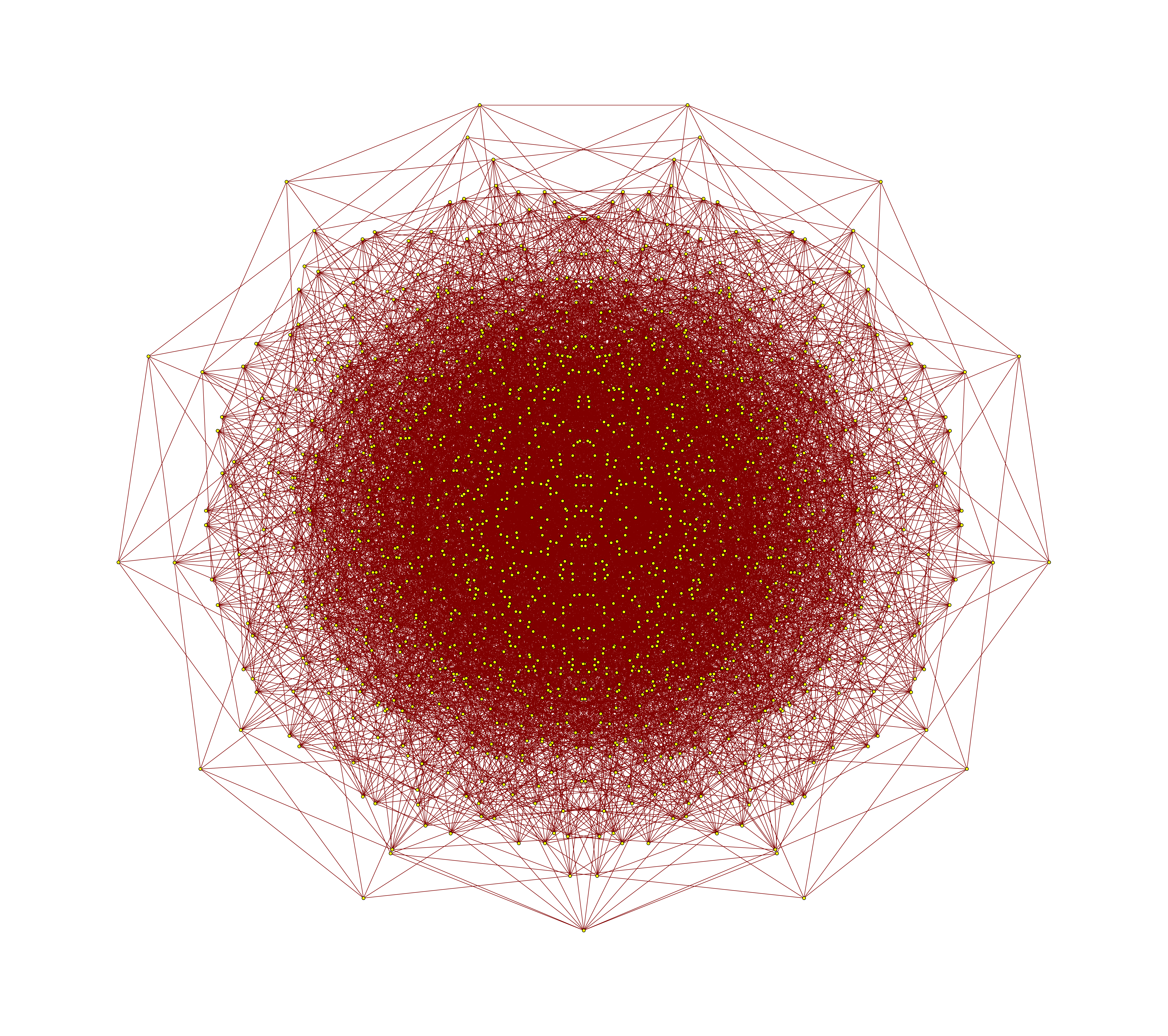}
\caption{The 1581-vertex, non-4-colourable unit-distance graph $G$.}
\label{ninthfig}
\end{figure}

\section{Acknowledgements}
I am indebted to Boris Alexeev, Rob Hochberg, Brendan McKay, Dustin Mixon, Paul Phillips, Landon Rabern and Gordon Royle for testing various graphs with independent code. I also thank Imre Leader and Graham Brightwell for comments on earlier versions of the manuscript, as well as for their role in germinating my interest in graph theory fully thirty years ago.

\bibliographystyle{amsplain}

\begin{thebibliography}{10}

\bibitem {Had} H. Hadwiger, \textit{Uberdeckung des euklidischen Raum durch kongruente Mengen}, Portugaliae Math. \textbf{4} (1945), 238--242.

\bibitem {MM} L. Moser and M. Moser, \textit{Solution to Problem 10},
Can. Math. Bull. \textbf{4} (1961), 187--189.

\bibitem {DHJP} D.H.J. Polymath, \textit{Polymath 16}, https://dustingmixon.wordpress.com/2018/04/14/polymath16-first-thread-simplifying-de-greys-graph/

\bibitem {Soi} A. Soifer, \textit{The Mathematical Coloring Book},
Springer, 2008, ISBN-13: 978-0387746401.

\end{thebibliography}

\end{document}